\documentclass{amsart}
\usepackage{amssymb, latexsym}
\theoremstyle{plain}
\newtheorem{theorem}{Theorem}
\newtheorem{corollary}{Corollary}
\newtheorem{lemma}{Lemma}

\theoremstyle{definition}

\begin{document}
\title[Generation of MUBs as powers of a unitary matrix]
{Generation of mutually unbiased bases \\ 
as powers of a unitary matrix in 2-power dimensions}
%\author{Konstantinos Drakakis${}^1$, Rod Gow${}^2$ and Scott Rickard${}^3$}

\author[Rod Gow]{Rod Gow}
\thanks{School of Mathematical Sciences,
University College,
Belfield, Dublin 4,
Ireland,
\emph{E-mail address:} \texttt{rod.gow@ucd.ie}}

\begin{abstract} Let $q$ be a power of 2. We show by representation theory that there exists a
$q\times q$ unitary matrix of multiplicative order $q+1$ whose powers generate $q+1$ complex pairwise mutually unbiased bases in $\mathbb{C}^q$.  When $q$ is a power of an odd prime, there is a $q\times q$ unitary matrix of multiplicative order $q+1$ whose first $(q+1)/2$ powers generate $(q+1)/2$ complex pairwise mutually unbiased bases. We also show how the existence of these matrices implies
the existence of special orthogonal decompositions of certain simple Lie algebras.
\end{abstract}
\maketitle
\section{Introduction}
\noindent 

\noindent  Let $d$ be a positive integer and let $\mathbb{C}^d$ denote a vector space of dimension $d$ over the field $\mathbb{C}$ of complex numbers. Let $\langle u,v \rangle $ denote a positive definite
hermitian form on  $\mathbb{C}^d\times \mathbb{C}^d$, which is linear in the first variable and 
conjugate linear in the second.  Let $B_1=\{\,u_1, \ldots, u_d\,\}$ and $B_2=\{\,v_1, \ldots, v_d\,\}$ be orthonormal bases of $\mathbb{C}^d$. We say that the bases $B_1$ and $B_2$ are \textit{mutually unbiased} if 
\[
|\langle u_i,v_j\rangle|^2=\frac{1}{d}
\]
for $1\le i,j\le n$. If we identify the orthonormal bases $B_1$ and $B_2$ with 
$d\times d$ unitary matrices $U_1$ and $U_2$, respectively, then the condition that the bases are mutually unbiased is equivalent to saying that each entry of the product $U_1U_2^\dagger=U_1U_2^{-1}$, where  $U_2^\dagger$ denotes the conjugate-transpose of $U_2$, has absolute value $1/\sqrt{d}$. 

Following \cite{KR}, let $N(d)$ denote the maximum number of orthonormal bases of $\mathbb{C}^d$ which are pairwise mutually unbiased. It has been proved that $N(d)\le d+1$, and furthermore  the
equality
$N(d)=d+1$ occurs whenever $d$ is a power of a prime. It is presently unknown if the equality $N(d)=d+1$ ever happens when $d$ is not a power of a prime. 

Explicit constructions of $d+1$ mutually unbiased bases when $d$ is a power of a prime have been based on properties of finite fields, \cite{KR}. Somewhat less explicit constructions use the irreducible representations of extra-special $p$-groups to realize the bases as eigenspaces of finite abelian groups acting on $\mathbb{C}^d$. See, for example, \cite{BBRV} and the references of \cite{KR}.

The main purpose of this paper is to show that when $d=2^a$, we can construct $d+1$ pairwise mutually unbiased bases of $\mathbb{C}^d$ as the powers of a unitary matrix of multiplicative order $d+1$. This matrix arises as an automorphism of an extra-special 2-group which has a special action on the maximal
abelian subgroups of the extra-special group. The equivalent construction when $d=p^a$, with $p$ an odd prime, yields a unitary matrix of multiplicative order $d+1$ whose first $(d+1)/2$ powers provide
$(d+1)/2$ pairwise mutually unbiased bases of $\mathbb{C}^d$. We show that it is generally impossible
to achieve $d+1$ pairwise mutually unbiased bases of $\mathbb{C}^d$ using the powers of a single matrix in the odd prime power case.
\medskip

\noindent
\textbf{Acknowledgement.}  We wish to thank Ferenc Sz\"oll\"osi of Budapest University of Technology and Economics for correspondence relating to questions raised in an earlier version of this paper.
In particular, Lemma 6 and Corollary 1 were inspired by his remarks.

 \section{Construction of a class of $p$-groups}
\noindent We shall try to keep our exposition reasonably self-contained and work from basic principles
of finite group theory.
Let $p$ be a prime and let $q=p^a$, where $a$ is a positive integer. Let $\mathbb{F}$ denote
the finite field of order $q^2$. We define a multiplication on the set $\mathbb{F}\times \mathbb{F}$
by putting
\[
(a,b)(c,d)=(a+c,a^qc+b+d)
\]
for all ordered pairs $(a,b)$ and $(c,d)$ in $\mathbb{F}\times \mathbb{F}$.
It is straightforward to see that $\mathbb{F}\times \mathbb{F}$ is a finite
group of order $q^4$, which we shall denote by
$G_q$. The identity element
is $(0,0)$ and the inverse of $(a,b)$ is $(-a,a^{q+1}-b)$.
The center $Z(G_q)$ of $G_q$ consists of all elements
$(0,b)$ and is elementary abelian of order $q^2$. 

We set $x=(a,b)$, $y=(c,d)$ and let $[x,y]$ denote the commutator $x^{-1}y^{-1}xy$. We find that
\[
[x,y]=(0, a^qc-c^qa).
\]
Note that any element $e$ of $\mathbb{F}$ expressible in the form $e=a^qc-c^qa$ satisfies $e^q=-e$.
It follows that the commutator subgroup $G_q'$, which is generated by the commutators $[x,y]$, consists of elements of the form $(0,e)$, where $e^q=-e$, and hence $|G_q'|=q$. 

Assuming that $a\ne 0$, it is easy to check that
the centralizer of $x$ consists of all elements $(c,d)$,
where $d$ is an arbitary element of $\mathbb{F}$ and $a^qc=ac^q$.  
Thus the centralizer of $x$ has order $q^3$. 

We now consider the irreducible complex characters of $G_q$. 

\begin{lemma} The group $G_q$ has exactly $q^3$ irreducible complex characters of degree $1$. All 
other irreducible complex characters of $G_q$ have degree $q$ and they vanish outside $Z(G_q)$.
\end{lemma} 

\begin{proof} 
Since the commutator subgroup
$G_q'$ has order $q$, $G_q$ has exactly $|G_q:G_q'|=q^3$ different irreducible characters of degree 1.
Now let $x$ be any non-central element of $G_q$. The second orthogonality relation implies that, as the centralizer of $x$ has order $q^3$,
\[
q^3=\sum |\chi(x)|^2,
\]
where the sum extends over all irreducible characters $\chi$ of $G_q$. The characters of degree 1 each contribute 1 to this sum, and since there are $q^3$ of them, we deduce that $\psi(x)=0$ for any irreducible character $\psi$ of degree greater than 1. Finally, Schur's Lemma implies that
$|\psi(z)|=\psi(1)$ for all elements $z$ of $Z(G_q)$ and hence we obtain by the first orthogonality relation
\[
|G_q|=q^4=\sum |\psi(z)|^2=|Z(G_q)|\psi(1)^2=q^2\psi(1)^2,
\]
where the sum extends over all elements $z$ of $Z(G_q)$. We deduce that $\psi(1)=q$, as required.
\end{proof}

Next, we construct an automorphism of $G_q$ which ultimately will be used to generate the mutually unbiased bases.

\begin{lemma} Let $\alpha$ be an element of order $q+1$ in $\mathbb{F}$. Then the mapping
$\sigma:G_q\to G_q$ given by
\[
\sigma(a,b)=(\alpha a,b)
\]
is an automorphism of order $q+1$ of $G_q$.
\end{lemma}

\begin{proof} Given elements $x=(a,b)$, $y=(c,d)$ in $G_q$, we have
\[
\sigma(xy)=\sigma(a+c,a^qc+b+d)=(\alpha(a+c), a^qc+b+d)
\]
and 
\[
\sigma(x)\sigma(y)=(\alpha a,b)(\alpha c,d)=(\alpha(a+c), \alpha^{q+1}a^qc+b+d).
\]
We deduce that $\sigma(xy)=\sigma(x)\sigma(y)$, since $\alpha^{q+1}=1$ by hypothesis. Thus,
$\sigma$ is an automorphism of $G_q$ and it is clear that $\sigma$ has order $q+1$, since
$\alpha$ has order $q+1$.
\end{proof}

Next, we construct a maximal abelian subgroup  $A$ of $G_q$ by setting
\[
A=\{(a,b): a^q=a\},
\]
where $b$ runs over $\mathbb{F}$.
It is straightforward to verify that $A$ is indeed abelian and $|A|=q^3$. Applying the automorphism
$\sigma$, we define subgroups $A_i$ so that
\[
A_i=\sigma^{i-1}(A)
\]
for $1\le i\le q+1$ (thus $A_1=A$). 

\begin{lemma}  Suppose that $p=2$. Then we have $A_i\cap A_j=Z(G_q)$ if $1\le i\ne j\le q+1$, and 
$G_q$ is the union of $A_1$, \dots, $A_{q+1}$.
\end{lemma}

\begin{proof} Suppose that $(c,d)\in A_i\cap A_j$, where $1\le i,j\le q+1$. Then we have $c=\alpha^{i-1}a=\alpha^{j-1} a_1$, where
$a^q=a$, $a_1^q=a_1$. Since we need only to investigate the case when  $a$ and $a_1$ are both
non-zero, we deduce that $(\alpha)^{q(i-j)}=\alpha^{i-j}$. However,  as $\alpha^q=\alpha^{-1}$,  we
obtain that $\alpha^{2(i-j)}=1$. Since $\alpha$ has odd order, this is only possible if $i=j$.
\end{proof}

The result above does not hold when $p$ is odd, since $\alpha^{(q+1)/2}=-1$, and consequently
\[
A_{i+(q+1)/2}=A_i
\]
for $1\le i\le (q+1)/2$. It is however easy to establish the following analogue of Lemma 3.

\begin{lemma}  Suppose that $p$ is odd. Then we have $A_i\cap A_j=Z(G_q)$ if 
$1\le  i\ne j\le (q+1)/2$.
\end{lemma}

\section{The irreducible representations of the group $G_q$}

\noindent Let $\chi$ be an irreducible character of $G_q$ of degree $q$ and let $X$ be an irreducible representation
of $G_q$ with character $\chi$. We may assume that $X(G_q)$ consists of unitary matrices. Lemma 1 shows that $\chi(x)=0$ if $x\not\in Z(G_q)$ and thus we see that no element outside $Z(G_q)$ is in the kernel of $X$. On the other hand, setting $Z=Z(G_q)$ for brevity,
$X(Z)$ is certainly contained in the center of $X(G_q)$ and it is a non-trivial cyclic group consisting of scalar multiples of the identity by Schur's Lemma. Since $Z$ is an elementary abelian $p$-group,
we deduce that $|X(Z)|=p$ and $|\ker X|=q^2/p$. It follows that $X(G_q)$ has order $q^2p$. Of course,
$X(G_q)$ is an extra-special $p$-group, but we do not need to know this.

Recall now the automorphism $\sigma$ of $G_q$. We define the conjugate representation 
$X^{\sigma}$ by 
\[
X^{\sigma}(x)=X(\sigma(x))
\]
for all $x\in G_q$. Since $\sigma$ fixes the elements of $Z$ pointwise, and the character of $X$ vanishes outside $Z$, it is clear that $X^{\sigma}$ also has character $\chi$. It follows that 
$X$ and $X^{\sigma}$ are equivalent representations and thus there exists an invertible matrix
$D$, say, satisfying
\[
X^{\sigma}(x)= X(\sigma(x))=D^{-1}X(x)D
\]
for all $x\in G_q$. We include a proof of the following well known lemma for completeness.

\begin{lemma} Assume the notation introduced above. Then multiplying the matrix $D$ by a
non-zero complex scalar if necessary, we may arrange for $D$ to satisfy $D^{q+1}=I$,
$\det D=1$,  and to be unitary.
Furthermore, the resulting $D$ has all its entries in the field of definition of $X$.
\end{lemma}

\begin{proof} We first note that we may assume that $D$ as defined above has its entries in the field
of definition of $X$, since it defines an equivalence between two representations over that
field. Now
it is straightforward to verify that
\[
D^{-i}X(x)D^i=X(\sigma^i(x))
\]
for each integer $i$ and hence 
\[
D^{-(q+1)}X(x)D^{q+1}=X(x)
\]
for all $x\in G_q$, since $\sigma$ has order $q+1$. 
Therefore,
$D^{q+1}$  commutes with all $X(x)$ and hence is a scalar
multiple of the identity by Schur's Lemma, say $D^{q+1}=\lambda I$. 

Let $d$ be the
determinant of $D$. Taking determinants in the equality $D^{q+1}=\lambda I$, 
we obtain 
\[
d^{q+1}=\lambda^q.
\]
We accordingly set $D_1=d\lambda^{-1}D$ and calculate that
\[
D_1^{q+1}=d^{q+1}\lambda^{-(q+1)}D^{q+1}=\lambda^{-1}D^{q+1}=I.
\]
Thus, if we replace $D$ by $D_1$, we obtain $D^{q+1}=I$. We note also that $\det D=1$ in this case.

Finally, since $D$ normalizes $X(G_q)$, and the elements $X(x)$ are unitary, we see that
\[
X(x)DD^\dagger X(x)^\dagger=X(x)DD^\dagger X(x)^{-1}=DD^\dagger
\]
Thus $DD^\dagger$ commutes with all elements $X(x)$ and is therefore a scalar multiple of the identity, by Schur's Lemma, say $DD^\dagger=\mu I$. But the diagonal entries of $DD^\dagger$
are sums of the squares of the moduli of complex numbers and are thus positive real numbers. This implies that $\mu$ is a positive real number. We also have $D^{q+1}=I$ and it follows that
\[
(D^\dagger)^{q+1}=I=\mu^{q+1}D^{-(q+1)}=\mu^{q+1}I.
\]
We deduce that $\mu^{q+1}=1$ and since $\mu$ is a positive scalar, $\mu=1$. This means that
$D$ is unitary.
\end{proof}

We note that we may take the field of definition of $X$ to be $\mathbb{Q}(\sqrt{-1})$ when $p=2$ and
$\mathbb{Q}(\omega)$, where $\omega$ is a primitive $p$-th root of unity, when $p$ is odd, and 
consequently we may assume that $D$ has all entries in the corresponding field.

\section{Construction of mutually unbiased bases using powers}

\noindent Recall that the abelian subgroups $A_i$ of $G_q$ satisfy $|A_i:Z(G_q)|=q$ for all $i$.
Let $x_{i,j}$, where $1\le j\le q$, be a system of coset representatives of $Z(G_q)$ in $A_i$, where
we choose $x_{i,1}$ to be the identity element.  Let $X$ be  a fixed irreducible complex representation of $G_q$ of degree $q$, consisting
of unitary matrices. 
We set
\[
\mathcal{C}_i=\{\,X(x_{i,j}):\ 1\le j\le q\,\}.
\]

The subsets $\mathcal{C}_i$ consist of $q$ commuting unitary matrices and satisfy $\mathcal{C}_i\cap \mathcal{C}_j=I$ for $i\ne j$ 

The trace inner product $f(U,V)$ of $q\times q$ complex matrices $U$ and $V$ is defined by
\[
f(U,V)=\hbox{tr} (U^\dagger V),
\]
where tr signifies the usual matrix trace. The trace inner product is a positive definite hermitian form
on the space of complex matrices. Since for any element $x$ of $G_q$ outside $Z(G_q)$ we have
$\hbox{tr}\,X(x)=0$, it is straightforward to see that the non-identity elements of $\mathcal{C}_i$ are
orthogonal to the non-identity elements of $\mathcal{C}_j$ with respect to the trace inner product.

We are now in a position to use a result of \cite{BBRV} to construct mutually unbiased bases.

\begin{theorem} Let $q$ be a power of $2$ and let $X$ be an irreducible complex representation of $G_q$ of degree $q$, consisting
of unitary matrices. Let $D$ be a $q\times q$ matrix that satisfies $D^{q+1}=I$ and  
$D^{-1}X(x)D=X(\sigma(x))$ for all $x$ in $G_q$. Then the powers $D$, $D^2$, \dots, $D^{q+1}=I$ define $q+1$ pairwise mutually unbiased bases. Furthermore, all entries of $D$ are in the field
$\mathbb{Q}(\sqrt{-1})$.
\end{theorem}

\begin{proof} The matrices in $\mathcal{C}_1$ commute and there therefore exists an orthonormal
basis of $\mathbb{C}^q$ consisting of simultaneous eigenvectors of the elements of $\mathcal{C}_1$.
This orthonormal basis determines a unitary matrix $B$, say. Now the matrix $D$ satisfies
\[
D^{-(i-1)}\mathcal{C}_1D^{i-1}=\mathcal{C}_{i}
\]
and thus the corresponding unitary matrix constructed of simultaneous orthonormal eigenvectors of
$\mathcal{C}_i$ is $D^{-(i-1)}B$. By Theorem 3.2 of \cite{BBRV}, the unitary matrices $D^{-(i-1)}B$ for
$1\le i\le q+1$ define $q+1$ pairwise mutually unbiased bases. But we may multiply these matrices on the right by $B^{-1}$ and still retain the unbiased property. The matrices so obtained are simply the powers of $D$.

\end{proof}
\medskip
\noindent\textbf{Example.} The matrix 
\[
\frac{1+i}{2}\left(
\begin{array} {rr}
-1&i\\
1& i
\end{array}
\right),
\]
where $i^2=-1$,
is unitary and has order 3. Its three different powers generate three mutually unbiased bases in
$\mathbb{C}^2$.

When $q$ is a power of an odd prime, we can obtain half the maximum number of mutually unbiased
bases from the powers of a unitary matrix, as a minor modification of the proof above shows.

\begin{theorem} Let $q$ be a power of an odd prime and let $X$ be an irreducible complex representation of $G_q$ of degree $q$, consisting
of unitary matrices. Let $D$ be a $q\times q$ matrix that satisfies $D^{q+1}=I$ and  
$D^{-1}X(x)D=X(\sigma(x))$ for all $x$ in $G_q$. Then the powers $I$, $D^{-1}$, \dots, $D^{-(q-1)/2}$ define $(q+1)/2$ mutually unbiased bases.
\end{theorem}

\begin{proof} The proof is identical with the previous one, except that, since
\[
D^{-(q+1)/2}\mathcal{C}_1D^{(q+1)/2}=\mathcal{C}_{1},
\]
we only obtain $(q+1)/2$ bases using the powers.
\end{proof}

There is an essential difference between the generation of mutually unbiased bases in the 2-power dimension case and the odd prime power dimension case, as we shall
soon show. We first prove a lemma on involutory
unitary matrices. 

\begin{lemma} Let $d$ be an odd positive integer and let $U$ be a $d\times d$ unitary matrix satisfying
$U^2=I$. Suppose that each diagonal entry of $U$ has absolute value $1/\sqrt{d}$. Then $d$ is a square.
\end{lemma} 

\begin{proof} Since we have
\[
UU^\dagger=U^2=I,
\]
it follows that $U=U^\dagger$, and thus $U$ is hermitian. The diagonal entries of $U$ are therefore real and consequently equal $\pm 1/\sqrt{d}$. As $U$ is an involution, its eigenvalues equal $\pm 1$ and hence its trace is an odd integer, since $d$ is odd by hypothesis. However the trace of $U$ equals the sum of its diagonal entries and thus we must have
\[
\frac{r-s}{\sqrt{d}}=e,
\]
where $e$ is an odd integer, $r$ is the number of diagonal entries of $U$ equal to $1/\sqrt{d}$, and
$s$ is the number of diagonal entries of $U$ equal to $-1/\sqrt{d}$. This implies that $\sqrt{d}$
is rational and hence $d$ is a square.
\end{proof}

\begin{corollary} Let $d$ be an odd positive integer which is not a square. Then there does not exist a $d\times d$ unitary matrix $D$ of even multiplicative order $2r$, say, whose powers $D$, $D^2$, \dots, $D^{2r}=I$ generate $2r$ pairwise mutually unbiased bases.
\end{corollary}

\begin{proof} Suppose on the contrary that such a unitary matrix $D$ exists. Then $D^r$ satifies the hypotheses of Lemma 6, which is impossible, since $d$ is not a square. Thus no such $D$ exists.
\end{proof}

We now see that there can be no analogue of Theorem 1 when $q$ is an odd power of an odd prime,
and thus we probably cannot improve Theorem 2.

\section{Reality questions for mutually unbiased bases}
\noindent There is a concept of real mutually unbiased bases, which appears to be much more restrictive than its complex counterpart. Suppose that we have a positive definite symmetric bilinear
form defined on $\mathbb{R}^d\times \mathbb{R}^d$. An orthonormal basis of $\mathbb{R}^d$
corresponds to a real $d\times d$ orthogonal matrix. Orthonormal bases $B_1$ and $B_2$ of $\mathbb{R}^d$ with 
corresponding $d\times d$ orthogonal matrices $O_1$ and $O_2$ are mutually unbiased when each entry of the product $O_1O_2^{-1}$ has absolute value $1/\sqrt{d}$. This means that
$\sqrt{d}O_1O_2^{-1}$ is an Hadamard matrix. This condition alone is enough to show that
pairs of real mutually unbiased bases can only exist for even values of $d$ (and, apart from $d=2$, only for values
of $d$ divisible by 4). See, for example, \cite{BSTW}, for information on the existence problem
for real mutually unbiased bases.

Let $N_{\mathbb{R}}(d)$ denote the maximum number of orthonormal bases of $\mathbb{R}^d$ which are pairwise mutually unbiased. It can be proved that $N_{\mathbb{R}}(d)\le d/2+1$, and furthermore  the
equality occurs when $d$ is a power of 4. As we hinted above, in many cases the value
of  $N_{\mathbb{R}}(d)$ for specific $d$ is substantially less than the general inequality
suggests. See \cite{BSTW} for much more information on this topic.

Constructions of families of real pairwise mutually unbiased bases meeting the upper bound
$d/2+1$ when $d$ is a power of 4 generally use orthogonal geometry over the field
of order 2 or real representations of extra-special 2-group. 
See, for example, Proposition 6 of \cite{CS} or Theorem 3.4 of \cite{CCKS}. We are interested in raising the question
of whether there exists any analogue of Theorem 1 for real mutually unbiased bases in this case.
Setting $q=2^{2m-1}$, we thus look for a real $2^{2m}\times 2^{2m}$ orthogonal matrix of order
$q+1$ whose powers generate $q+1$ real mutually unbiased bases. At present, we have not succeeded in showing that such matrices exist, but we note the following $4\times 4$ example.
\medskip

\noindent\textbf{Example.} The matrix 
\[
\frac{1}{2}\left(
\begin{array} {rrrr}
-1&-1&-1&-1\\
1&-1&-1&1\\
1&1&-1&-1\\
1& -1 &1&-1
\end{array}
\right)
\]
is orthogonal and has order 3. Its three different powers generate three mutually unbiased bases in
$\mathbb{R}^2$. We obtain an Hadamard matrix if the matrix above is multiplied by 2.

\section{Connections with orthogonal decompositions of simple Lie algebras}

\noindent The paper \cite{BSTiepW} has shown a connection between the existence of $d+1$ mutually unbiased bases in $\mathbb{C}^d$ and orthogonal decompositions of the special
linear Lie algebra $sl_d(\mathbb{C})$ into a direct sum of Cartan subalgebras. In this section we show
how the automorphism $\sigma$ introduced above gives extra information on such a decomposition
in the 2-power case. We also obtain similar information about an orthogonal decomposition of the symplectic Lie algebra of $d\times d$ matrices when $d$ is a 2-power.

We employ the notation of Section 4. Let $\mathcal{L}$ denote the special linear Lie algebra
of all $q\times q$ complex matrices of trace zero. Recall that we have defined the subsets $\mathcal{C}_i$, which consist of $q-1$ commuting matrices of trace zero, together with the identity. These trace zero matrices are linearly independent, since they are pairwise orthogonal with respect to the trace inner product.
Let $\mathcal{H}_i$ denote the linear span of the non-identity elements in $\mathcal{C}_i$. Clearly,
$\mathcal{H}_i$ is an abelian subalgebra of $\mathcal{L}$ of dimension $q-1$. Furthermore, as indicated in \cite{BSTiepW}, Theorem 5.1, $\mathcal{H}_i$ is self-normalizing in $\mathcal{L}$. For,
up to conjugation by a unitary matrix, which induces an automorphism of $\mathcal{L}$, we may assume
that $\mathcal{H}_i$ is a $q-1$-dimensional subspace of diagonal trace zero matrices, and hence 
equals the entire subspace of all such diagonal matrices, which is  well known to be a Cartan subalgebra. 
We note also that $\mathcal{H}_i$ is orthogonal to $\mathcal{H}_j$ with respect to the Killing form of
$\mathcal{L}$, since this is a scalar multiple of the usual trace form, \cite{BSTiepW}, Lemma 4.2.

We therefore  obtain a decomposition
\[
\mathcal{L}=\mathcal{H}_1\perp\cdots \perp\mathcal{H}_{q+1}
\]
of $\mathcal{L}$ into an orthogonal direct sum of $q+1$ Cartan subalgebras. Now conjugation
by the matrix $D$ described in Theorem 1 induces an automorphism of order $q+1$ of
$\mathcal{L}$. In the case that $q$ is a power of 2, this automorphism transitively permutes
the $q+1$ Cartan subalgebras, by the arguments of Theorem 1. (The automorphism of course preserves the Killing form.) Thus we may summarize our findings
in this case.

\begin{theorem} Let $q$ be a power of $2$ and let $\mathcal{L}$ denote the Lie algebra of all
$q\times q$ complex matrices of trace $0$. Then there is an orthogonal decomposition
\[
\mathcal{L}=\mathcal{H}_1\perp\cdots \perp\mathcal{H}_{q+1}
\]
of $\mathcal{L}$ into a direct sum of $q+1$ Cartan subalgebras $\mathcal{H}_i$ which are transitively
permuted by an automorphism of $\mathcal{L}$ of order $q+1$.
\end{theorem}

We observe that a larger finite group of automorphisms
permutes the summands of the decomposition above, while preserving the Killing form, since the group generated by $X(G_q)$ and $D$ also has this property.

This result bears a superficial resemblance to a deep result of Thompson, which shows that the
complex Lie algebra $\mathcal{L}$ of type $E_8$ admits an orthogonal decomposition into a direct
sum of 31 Cartan subalgebras, which are transitively permuted by a group of automorphisms
isomorphic to the Dempwolff group. 

We do not know if the corresponding orthogonal decomposition of the special linear Lie algebra
admits a similar transitive action of an automorphism of order $q+1$ when $q$ is a power of an odd prime.

The decomposition described in Theorem 3 can be refined into a similar decomposition of the
symplectic Lie algebra. We sketch some of the details here. The representation $X$ of $G_q$ described in Theorem 1 has a real-valued character but it cannot be defined over the real numbers. We say that it is of symplectic type, as the elements $X(x)$ preserve a non-degenerate symplectic form. This means that
there is an invertible $q\times q$ skew-symmetric matrix $S$, say, such that 
\[
X(x)SX(x)^T=S
\]
for all $x\in G_q$. Furthermore, it is straightforward to see that, given our normalization of $D$,
$D$ also satisfies $DSD^T=S$.  

The symplectic Lie algebra $\mathcal{L}$ consists of all $q\times q$ matrices $A$ which satisfy
\[
AS+SA^T=0.
\]
It has dimension $q(q+1)/2$ and rank $q/2$.

Suppose that $x$ is an element of $G_q$ such that $X(x)$ has order 4. Then $X(x)^2=-I$ and hence
\[
X(x)^{-1}=-X(x).
\]
It follows that the equality $X(x)SX(x)^T=S$ is equivalent to $X(x)S+SX(x)^T=0$ in this case. Thus the
elements of order 4 in $X(G_q)$ are in the symplectic Lie algebra. 

Consider now the abelian subgroup
$X(A_i)$ of $X(G_q)$. It is straightforward to show that, if $q=2^m$, it is isomorphic to a direct product
\[
\mathbb{Z}_4\times \mathbb{Z}_2\times\cdots\times \mathbb{Z}_2
\]
of cyclic groups and thus contains $2^m$ elements of order 4.  Let $\mathcal{D}_i$ denote the linear
span of the elements of order 4 in $X(A_i)$. This is an abelian subalgebra of dimension $q/2$ of  $\mathcal{L}$, which must be a Cartan subalgebra, as it consists of diagonalizable
elements and has the appropriate dimension. 

The $q+1$ subalgebras  $\mathcal{D}_i$, where
$1\le i \le q+1$,  are all orthogonal with respect to the Killing form of $\mathcal{L}$, for the same reason as given in the proof of Theorem 3. Finally, $D$ induces an automorphism
of order $q+1$ of $\mathcal{L}$ by conjugation, since it preserves the symplectic form determined by
$S$ (indeed the group generated by $X(G_q)$ and $D$ has this property). This automorphism
transitively permutes the Cartan subalgebras $\mathcal{D}_i$, since  $D$
maps the elements of order 4 in $X(A_i)$ into the elements of order 4 in
$X(A_{i+1})$. 

Summarizing, we have proved the following.

\begin{theorem} Let $q$ be a power of $2$ and let $\mathcal{L}$ denote the 
symplectic Lie algebra of 
$q\times q$ complex matrices. Then there is an orthogonal decomposition
\[
\mathcal{L}=\mathcal{H}_1\perp\cdots \perp\mathcal{H}_{q+1}
\]
of $\mathcal{L}$ into a direct sum of $q+1$ Cartan subalgebras $\mathcal{H}_i$ which are transitively
permuted by an automorphism of $\mathcal{L}$ of order $q+1$.
\end{theorem}

It is interesting to observe that at the time of this writing, orthogonal decompositions of the symplectic
Lie algebra are only known in this 2-power rank case. The paper \cite{K} contains related and more
general constructions for orthogonal decompositions of Lie algebras using the geometry of finite
vector spaces,
and contains many references to earlier work.

\end{document}